\documentstyle[12pt]{article}

\reversemarginpar
\newcommand{\nwc}{\newcommand}
\newcommand{\commentout}[1]{}
\newcommand{\bes}{\begin{displaymath}}
\newcommand{\ees}{\end{displaymath}}
\newcommand{\be}{\begin{equation}}
\newcommand{\ee}{\end{equation}}
\newcommand{\ba}{\begin{eqnarray}}
\newcommand{\ea}{\end{eqnarray}}
\newcommand{\bas}{\begin{eqnarray*}}
\newcommand{\eas}{\end{eqnarray*}}

\newcommand{\al}{\alpha}
\newcommand{\bbs}{{\bf s}}
\newcommand{\bbsig}{{\bf \si}}
\newcommand{\bbj}{{\bf j}}
\newcommand{\bE}{{\bf E}}
\newcommand{\bbl}{{\bf l}}
\newcommand{\cF}{{\cal F}}
\newcommand{\cG}{{\cal G}}

\newcommand{\bI}{{\bf I}}
\newcommand{\bW}{{\bf W}}

\newcommand{\bD}{{\bf D}}
\newcommand{\bone}{{\bf 1}}

\newcommand{\bbn}{{\bf n}}
\newcommand{\bt}{\beta}

\newcommand{\si}{\sigma}

\newcommand{\bB}{{\bf B}}

\newcommand{\vV}{{\bf V}}

\newcommand{\vphi}{\varphi}

\newcommand{\ep}{\varepsilon}

\newcommand{\Om}{\Omega}

\newcommand{\bR}{{\bf R}}

\newcommand{\bk}{{\bf k}}

\newcommand{\bze}{{\bf 0}}

\newcommand{\kropa}{$_{\Box}$}

\newcommand{\bbi}{{\bf i}}
\newcommand{\bbt}{{\bf t}}

\newcommand{\da}{\downarrow0}
\newtheorem{theorem}{Theorem}
\newtheorem{proposition}{Proposition}

\newtheorem{lemma}{Lemma}

\nwc{\m}{\mbox}
\nwc{\bm}{\boldmath}
\nwc{\ubm}{\unboldmath}
\nwc{\bmu}{\m{\bm $u$\ubm}}
\nwc{\bmx}{\m{\bm $x$\ubm}}
\nwc{\bx}{\bmx}
\newcommand{\by}{\m{\bm $y$\ubm}}
\nwc{\bmv}{\m{\bm $v$\ubm}}
\nwc{\cE}{{\cal E}}

\nwc{\beq}{\begin{eqnarray}}
\nwc{\eeq}{\end{eqnarray}}

\begin{document}
\title{Fractional Brownian Motion Limit for a Model
of Turbulent Transport}\author{Albert Fannjiang \\ Department of
Mathematics, University of California, Davis\and Tomasz Komorowski \\ 
Institute of Mathematics,
Maria Curie-Sk\l odowska University, Lublin
\thanks{The research of A. F. is supported by National Science Foundation
Grant No. DMS-9600119.}
}

\maketitle

{\bf Abstract}  Passive scalar motion in a family of
random Gaussian velocity fields
 with long-range correlations  
 is shown to converge to persistent fractional Brownian motions
 in long times.



 {\bf Abbreviated title} fractional Brownian motion limit.

\newpage
%\tableofcontents
\section{Introduction}
We consider the motion of a passive scalar advected by a random velocity
field $\vV(t,\bmx)=(V_1(t,\bmx),\cdots,V_d(t,\bmx))$.
The governing equation is
\be
\label{a1.0}
\frac{d\bmx(t)}{dt}=\vV(t,\bmx(t))
\ee
where $\vV(t,\bmx)$ is a mean-zero, time-stationary, space-homogeneous
random incompressible velocity field.

In certain situations, it is believed that the {\em convergence}
of the Taylor-Kubo formula (\cite{T}, \cite{K}) given by
\be
\label{TK}
\int\limits^\infty_0\left\{\bE[V_i(t,\bze)V_j(0,\bze)]+
\bE[V_i(t,\bze)V_j(0,\bze)]\right\}dt
\ee
is a criterion for convergence of passive scalar motion
to Brownian motion in the long time limit. Indeed, it has been
shown that the solution of
\be
{d\bmx^\ep(t)\over dt}={1\over\ep}\vV({t\over\ep^2},\bmx^\ep(t)),
\quad  \bmx^\ep(0)=0
\label{BM}
\ee
converges in law, as $\ep\to 0$,
to the Brownian motion with diffusion coefficients given by
the Taylor-Kubo formula when the velocity field is sufficiently {\em mixing
in time} (see \cite{khasminskii}, \cite{kesten-papanicolaou},
\cite{kunita}, \cite{carmona-fouque}).
Moreover, the solution of (\ref{BM})
converges to
the same Brownian motion for a family of {\em non-mixing}
Gaussian, Markovian flows with power-law spectra
as long as the Taylor-Kubo formula converges
(see \cite{fannjiang-komorowski}).  In this paper,
for the same family
of power-law spectra, we show that,
when the Taylor-Kubo formula {\em diverges},
the solution of the following equation
\be
\label{a1.1}
{d\bmx_\ep(t)\over dt}= \ep^{1-2\delta}\vV({t\over \ep^{2\delta}},\bmx_\ep(t)),
\quad \bmx^\ep(0)=0,
\ee
with some $\delta\neq 1$ depending on the velocity spectrum,
converges, as $\ep \to 0$,
to a  fractional Brownian motion (FBM), as introduced
in \cite{MV} (see also \cite{ST}).

We define the family of
velocity fields with power-law spectra as follows.
Let $(\Om,{\cal V},P)$ be a probability space of which
each element is a velocity field
$\vV(t,\bx)$, $(t,\bx)\in R\times R^d$
satisfying the following properties.
\begin{itemize}
\item[H 1)] $\vV(t,\bx)$ is time stationary, space-homogeneous and
centered, i.e., $\bE\{\vV\}=\bze$, and Gaussian.
Here $\bE$ stands for the expectation
with respect to the probablity measure $P$.
\item[H 2)] The two-point correlation tensor $\bR=[R_{ij}]$ is given by
\be
\label{a1.3}
R_{ij}(t,\bmx)=\bE\left[V_i(t,\bmx)V_j(0,\bze)\right]=
\int_{R^d}\cos{(\bk\cdot\bx)}e^{-|\bk|^{2\bt}t}\hat{\bR}_{ij}(\bk)d\bk
\ee
with the spatial spectral density
\be
\label{060103}
\label{a1.3b}
\hat{\bR}(\bk)=\frac{a(|\bk|)}{|\bk|^{2\al+d-2}}\left(\bI-\frac{\bk\otimes\bk}{|
\bk|^2}\right),
\ee
where $a:[0,+\infty)\rightarrow R_+$ is a compactly supported, continuous,
nonnegative function.
The factor $\bI-\bk\otimes\bk/|\bk|^2$ in (\ref{a1.3b}) is a result of
incompressibility.

\item[H 3)] $\al<1$, $\bt\geq 0$ and $\al+\bt>1$.
\end{itemize}

The function $\exp{(-|\bk|^{2\beta}t)}$ in (\ref{a1.3}) is called the
{\em time correlation function} of the flow $\vV$.
For $\beta>0$, the velocity field
lacks the spectral gap and, thus, is not mixing in time.
As the time correlation function is exponential,
the Gaussian velocity field is Markovian in time.

Because the function $a$ has a compact
support we may assume, without loss of generality, that $\vV$ is jointly
continuous in both $(t,\bx)$ and is $C^\infty$ in $\bx$ almost surely.
For $\alpha<1$, the spectral density $\hat{\bR}(\bk)$ is integrable in $\bk$
and, thus, (\ref{a1.3})-(\ref{060103}) defines a random
velocity field with a finite
second moment. The exponent $\alpha$ is directly related to the decay exponent
of $\bR$. Namely $|\bR|(0,\bmx)\sim |\bmx|^{\al-1}$ for $|\bmx|\gg1$. As
$\alpha$ increases to one, the decay exponent of $\bR$ decreases to zero.

\begin{figure}
\begin{picture}(400,250)
\includegraphics{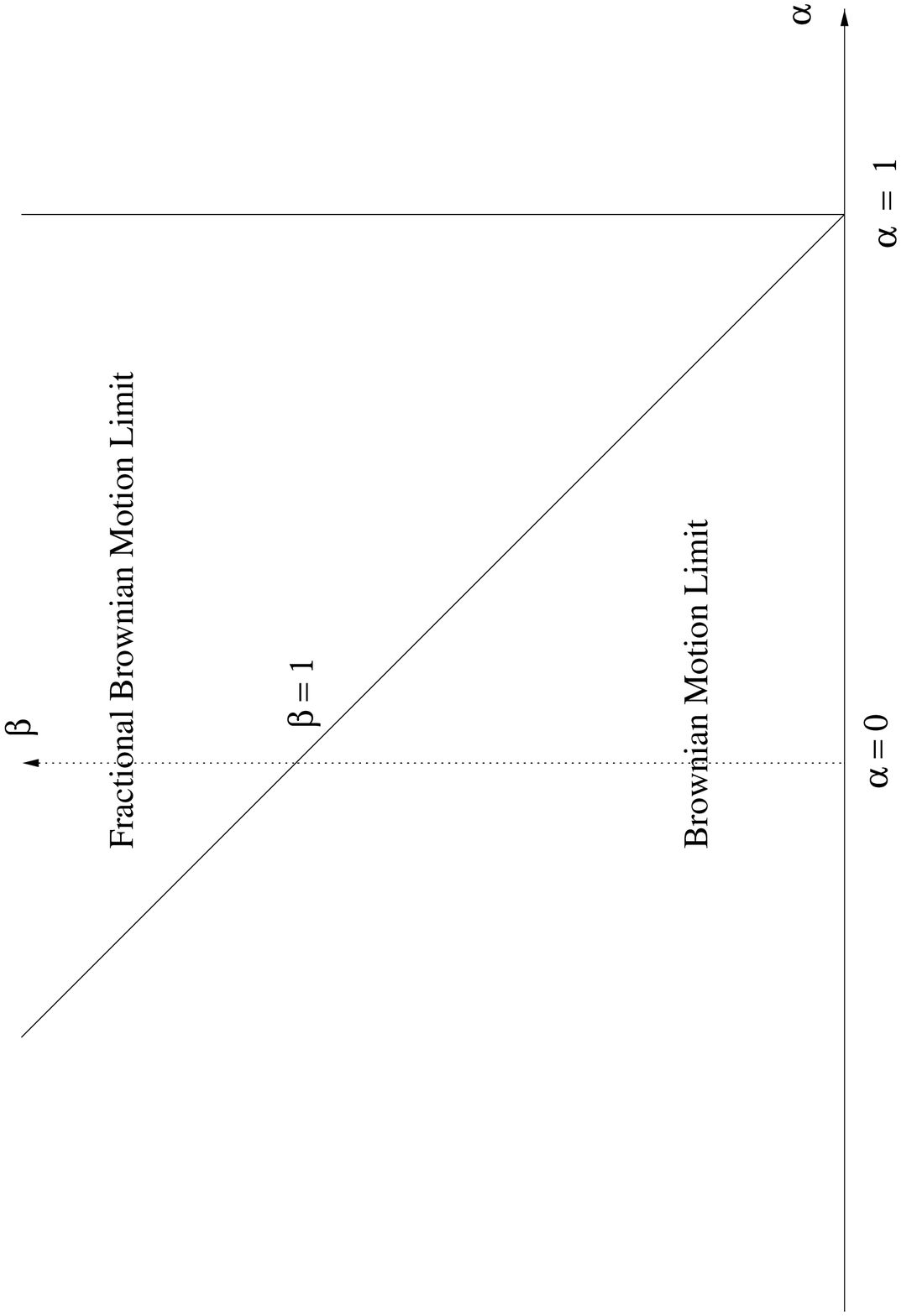}
\end{picture}
\end{figure}

Our main result is summarized in the
 following theorem.
\begin{theorem}
\label{theorem1}
Under the assumptions H 1)- H 3),
the solution of eq. (\ref{a1.1}) with the scaling exponent
\[
\delta:=\frac{\bt}{\al+2\bt-1}
\]
converges in law, as $\ep$ tends to zero, to a fractional Brownian motion
$\bB_H(t)$ that is to a Gaussian process with stationary increments whoe
covariance is given by
\be
\label{062701}
\bE\left[\bB_H(t)\otimes\bB_H(t)\right]=\bD t^{2H},
\ee
with the coefficients $\bD$
\be
\label{061802}
\bD=\int\limits_{R^d}\frac{
e^{-|\bk|^{2\bt}}-1+|\bk|^{2\bt}}{|\bk|^{2\al+4\bt-1}}
\left(\bI-\frac{\bk\otimes\bk}{|\bk|^2}\right){a(0)\over |\bk|^{d-1}}d\bk
\ee
and the Hurst exponent $H$
\be
\label{061803}
1/2<H=1/2+\frac{\al+\bt-1}{2\bt}<1.
\ee
\end{theorem}
{\bf Remark.} Molecular diffusion can be added to the equation of motion
 so that
instead of (\ref{a1.0}) we may consider an It\^{o} stochastic differential
equation
\[
d\bmx(t)=\vV(t,\bmx(t))dt+\sqrt{2\kappa}d\bB(t)
\]
with $\bB(t)$, $t\geq0$ the standard Brownian motion, independent of $\vV$ and
$\kappa\geq0$. This however would not influence our results.

\section{Multiple stochastic integrals}

By the Spectral Theorem (see,  e.g., \cite{adler}) we assume
without loss of any generality that there exist two
independent, identically distributed, real vector valued, Gaussian spectral measures
$\hat{\vV}_{l}(t,\cdot)$, $l=0,1$ such that
\be
\vV(t,\bx)=\int \hat{\vV}_0(t,\bx,d\bk),
\label{a.0}
\ee
where
\[
\hat{\vV}_0(t,\bx,d\bk):=c_0(\bk\cdot\bx)\hat{\vV}_0(t,d\bk)+c_1(\bk\cdot\bx)\hat{\vV}_1(t,d\bk)
\]
with $c_0(\phi)\equiv\cos{(\phi)},c_1(\phi)\equiv\sin{(\phi)}$.
Define also
\[
\hat{\vV}_1(t,\bx,d\bk):=-c_1(\bk\cdot\bx)\hat{\vV}_0(t,d\bk)+c_0(\bk\cdot\bx)\hat{\vV}_1(t,d\bk).
\]
We have the relations
\beq
\label{a.2}
\partial \hat{\vV}_0(t,\bx,d\bk)/\partial x_j&= &k_j \hat{\vV}_1(t,\bx,d\bk),\\
\quad\partial \hat{\vV}_1(t,\bx,d\bk)/\partial x_j&=& -k_j \hat{\vV}_0(t,\bx,d\bk).\label{a.1}
\eeq
Clearly $\int \hat{\vV}_1(t,\bx,d\bk)$
is a random field distributed identically to and independently of $\vV$.

We define the {\em multiple stochastic integral}
\be
\label{08144}
\int\cdots\int  \psi(\bk_1,\cdots,\bk_N)\widehat{\vV}_{l_1}(t_1,\bx_1,d\bk_1)\otimes\cdots\otimes\widehat{\vV}_{l_N}(t_N,\bx_N,d\bk_N)
\ee
for any  $l_1,\cdots,l_N\in\{0,1\}$ and a suitable family of functions $\psi$
by using  the Fubini theorem (see (\ref{061902}) below).
For $\psi_1,\cdots,\psi_N\in {\cal S}(R^d)$,
the Schwartz space, and $l_1,\cdots,l_N\in\{0,1\}$  we set
\beq
\label{061902}
&&
\int\cdots\int \psi_1
(\bk_1)\cdots\psi_N(\bk_N)\widehat{\vV}_{l_1}(t_1,\bx_1,d\bk_1)\otimes\cdots\otimes\widehat{\vV}_{l_N}(t_N,\bx_N,d\bk_N)\\
&:=&
\int\psi_{1}(\bk_1)\widehat{\vV}_{l_1}(t_1,\bx_1,d\bk_1)\otimes\cdots\otimes\int\psi_N(\bk_N)\widehat{\vV}_{l_N}(t_N,\bx_N,d\bk_N).\nonumber
\eeq
We then extend the definition of multiple integration to the closure ${\cal H}$ of the Schwartz space ${\cal S}((R^d)^N,R)$ under the norm
\be
\label{0807}
\|\psi\|^2:=\int\cdots\int \psi(\bk_1,\cdots,\bk_N) \psi(\bk_1',\cdots,\bk_N')
\ee
\[
\bE\left[\widehat{\vV}_{l_1}(t_1,\bx_1,d\bk_1)\otimes\cdots\otimes
\widehat{\vV}_{l_N}(t_N,\bx_N,d\bk_N)\cdot\widehat{\vV}_{l_1}(t_1,\bx_1,d\bk_1')\otimes\cdots\otimes\widehat{\vV}_{l_N}(t_N,\bx_N,d\bk_N')\right].
\]
The expectation is to be calculated by the formal rule
\[
\bE\left[\widehat{V}_{l,i}(t,\bx,d\bk)\widehat{V}_{l',i'}(t',\bx',d\bk')\right]=e^{-|\bk|^{2\bt}|t-t'|}\delta_{l,l'}c_0(\bk\cdot(\bx-\bx'))
\widehat{R}_{i,i'}(\bk)\delta(\bk-\bk')d\bk d\bk'.
\]
This approach to spectral integration follows \cite{shiryaev}.

When  $\bbi=(i_1,\cdots,i_d),i_1,\cdots,i_d\in
\{1,2,\cdots,d\}$ is fixed and
$\bbl=(l_1,\cdots,l_N), l_1,\cdots,l_N\in
\{0,1\}$ we shall denote the corresponding component of
the stochastic integral by $\Psi_{\bbl,\bbi}$.

Note that $\Psi_{\bbl,\bbi}\in
H^N(\vV)$ - the Hilbert space
obtained as a completion of the space of $N$-th degree polynomials in
variables $\int\psi(\bk)\widehat{\vV}(t,\bx,\bk)$ with
respect to the standard $L^2$ norm.
\begin{proposition}
\label{proposition611}
For any $(t_1,\bx_1),\cdots,(t_N,\bx_N)\in R\times R^d$ and $p>0$,
$\Psi_{\bbl,\bbi}$ belongs to  $L^p(\Om)$ and
\be
\label{061701}
\left(\bE|\Psi_{\bbl,\bbi} |^p\right)^{1/p}\leq C\left(\bE|\Psi_{\bbl,\bbi}|^2\right)^{1/2}
\ee
with the constant $C$ depending only on $p,N$ and the dimension $d$.
Moreover, $\Psi_{\bbl,\bbi}$ is differentiable in the mean square sense with
\be
\label{061702}
\nabla\Psi_{\bbl,\bbi}(t_1,\cdots,t_N,\bx_1,\cdots,\bx_N)=
(-1)^{l_j}\int\cdots\int \bk_j \psi(\bk_1,\cdots,\bk_N)
\ee
\[
\widehat{V}_{l_1,i_1}(t_1,\bx_1,d\bk_1)\cdots\widehat{V}_{1-l_j,i_j}(t_j,\bx_j,d\bk_j)\cdots\widehat{V}_{l_N,i_N}(t_N,\bx_N,d\bk_N).
\]
\end{proposition}

The proof of Proposition
\ref{proposition611} is standard and follows directly from the well known
hypercontractivity property for Gaussian measures
(see, e.g., \cite{janson}, Theorem 5.1. and its corollaries), so
we do not repeat it here.

The field $\vV$ is Markovian i.e.
\be
\label{1.10}
\bE\left[\int\psi(\bk)\widehat{\vV}_l(t,\bx,d\bk)\left|\right.{\cal
    V}_{-\infty,s}\right]= \int
e^{-|\bk|^{2\bt}(t-s)}\psi(\bk)\widehat{\vV}_l(s,\bx,d\bk),~~l=0,1,
\ee
 for  all
 $\psi\in{\cal S}(R^d,R)$,
where ${\cal
    V}_{a,b}$ denotes the $\si$-algebra generated by random variables
    $\vV(t,\bx)$, for $t\in [a,b]$ and $\bx\in R^d$.

To calculate a mathematical expectation of multiple product of Gaussian random
variables, it is convenient to use a graphical representation, borrowed from
quantum field theory. We refer to, e.g., Glimm and Jaffe \cite{GJ} and Janson
\cite{janson}. A {\em Feynman diagram} ${\cal F}$ (of order $n\geq 0$ and rank
$r\geq0$) is a graph consisting of a set $B(\cF)$ of $n$ vertices and a set
$E(\cF)$ of $r$ edges without common endpoints. So there are $r$ pairs of
vertices, each joined by an edge, and $n-2r$ unpaired vertices, called  {\em
free vertices}. $B(\cF)$ is a set of positive integers. An edge whose endpoints
are $m,n\in B$ is represented by $\widehat{mn}$ (unless otherwise specified, we
always assume $m<n$); and an edge includes its endpoints. A diagram ${\cal F}$
is said to be {\em based on} $B(\cF)$. Denote the set of free vertices  by
$A({\cal F})$, so $A({\cal F})=\cF\setminus E(\cF)$. The diagram is {\em
complete} if $A({\cal F})$ is empty and {\em incomplete}, otherwise. Denote by
$\cG(B)$ the set of all diagrams based on $B$, by $\cG_c(B)$ the set of all
complete diagrams based on $B$ and by $\cG_i(B)$ the set of all incomplete
diagrams based on $B$. A diagram ${\cal F}'\in \cG_c(B)$ is called a {\em
completion} of ${\cal F}\in\cG_i(B)$ if $E({\cal F})\subseteq E({\cal F}')$.

Let $B=\{1,2,3,...,n\}$. Denote by ${\cal F}_{|k}$ the sub-diagram of ${\cF}$,
based on $\{1,\cdots,k\}$. Define $A_k({\cal F})=A({\cal F}_{|k})$. A special
class of  diagrams, denoted by $\cG_s(B)$, plays an important role  in the
subsequent analysis: a diagram $\cF$ of order $n$ belongs to $\cG_s(B)$ if
$A_k(\cF)$ is not empty for all $k=1,...,n.$

We shall adopt the following multiindex notation. For any $P\in Z^+$, multiindex
$\bbn=(n_1,\cdots,n_P)$, $|\bbn|$ stands for $\sum n_p$. If  $P'\leq P$ we
denote $\bbn_{|P'}:=(n_1,\cdots,n_{P'})$. In addition if $k$ is any number we
set $\bbn\cdot k:=(n_1,\cdots,n_P,k)$.

 We work out the conditional expectation for
multiple spectral integrals using the Markov property (\ref{1.10}).
\begin{proposition}
For any function $\psi\in {\cal H}$ and $l_1,\cdots,l_N\in \{0,1\}$,
$i_1,\cdots,i_N\in\{1,\cdots,d\}$,
\label{prop1}
\be
\label{060106}
\bE\left[\int\cdots\int
\psi(\bk_1,\cdots,\bk_N)
\widehat{V}_{l_1,i_1}(t,\bx_1,d\bk_1)\cdots
\widehat{V}_{l_N,i_N}(t,\bx_N,d\bk_N)\left|\right.{\cal V}_{-\infty,s}\right]=
\ee
\[
\sum\limits_{{\cal
   F}\in\cG(\{1,...,N\})}\int\cdots\int\exp\left\{-\sum\limits_{m\in A({\cal
   F})}|\bk_m|^{2\bt}(t-s)\right\}\psi(\bk_1,\cdots,\bk_N)
\widehat{V}_{s,\bx_1,\cdots,\bx_N}(d\bk_1,\cdots,d\bk_N;{\cal
    F})
\]
with
\be
\label{01201}
\widehat{V}_{s,\bx_1,\cdots,\bx_N}(d\bk_1,\cdots,d\bk_N;{\cal
   F}):=\prod\limits_{m\in A({\cal
  F})}\widehat{V}_{l_m,i_m}(s,\bx_m,d\bk_m)
\ee
\[
\times\mathop{\prod\limits_{\widehat{mn}\in
   E({\cal F})}}
\left[1-e^{-\left(|\bk_{m}|^{2\bt}+|\bk_{n}|^{2\bt}\right)(t-s)}
\right]\bE[\widehat{V}_{l_{m},i_m}(s,\bx_m,d\bk_{m})
\widehat{V}_{l_{n},i_n}(s,\bx_n,d\bk_{n})].
\]
\end{proposition}

{\bf Proof.} Without loss of generality we consider
$
\psi(\bk_1,\cdots,\bk_N)=\bone_{A_1}(\bk_1)\cdots\bone_{A_N}(\bk_N)
$
for some Borel sets $A_1,\cdots,A_N$.

Note that
$\widehat{\vV}_l(t,A_i)=\widehat{\vV}^{0}_l(t,A_i)+\widehat{\vV}^{1}_l(t,A_i)$ where
$\widehat{\vV}^{0}_l(t,\cdot)$ is
the orthogonal projection of $\widehat{\vV}_l(t,\cdot)$ on
$L^2_{-\infty,t}$ and $\widehat{\vV}^{1}_l(t,\cdot)$ its complement.
Here  $L^2_{a,b}$ denotes $L^2$
closure of the linear span over $\vV(s,\bx)$, $a\leq s\leq b$, $\bx\in R^d$.
The conditional expectation in
(\ref{060106}) equals
\[
\sum\limits_{{\cal
    F}\in\cG(\{1,...,N\})}\mathop{\prod\limits_{\widehat{mn}\in
       E({\cal F})}}
     \bE\left[\widehat{V}_{i_{m},l_{m}}^1(t,A_{m})\widehat{V}_{i_{n},l_{n}}^1(t,A_{n})\right]\prod\limits_{m\in A({\cal F})}\widehat{V}_{i_m,l_m}^0(t,A_{m}).
\]
The statement follows upon the application of the relations
\[
\widehat{\vV}_{l}^0(t,A)=\int\limits_{A}e^{-|\bk|^{2\bt}(t-s)}\widehat{\vV}_{l}(s,d\bk)
\]
and
\[
\bE\left[\widehat{\vV}_{l}^1(t,A)\otimes\widehat{\vV}_{l'}^1(t,B)\right]=
\]
\[
\int\limits_{A}\int\limits_{B}\delta_{l,l'}\left\{\bE\left[\widehat{\vV}_{l}(t,d\bk)\otimes\widehat{\vV}_{l'}(t,d\bk')\right]-\bE\left[\widehat{\vV}_{l}^0(t,d\bk)\otimes\widehat{\vV}_{l'}^0(t,d\bk')\right]\right\}.
\]
\quad\rule{2mm}{3mm}

\label{sec2}

\section{Proof of tightness.}
\label{section2}
We begin with the following lemma
which shows, among other things, that the family of continuous trajectory
processes $\bx_\ep(t)$, $t\geq0$ is tight.
\begin{lemma}
\label{lemmab1}
For the family of trajectories given by (\ref{a1.1}) we have
\[
\lim\limits_{\ep\da}\bE\left[(\bx_\ep(t)-\bx_\ep(\tau))\otimes(\bx_\ep(t)-\bx_\ep(\tau))\right]=\bD
|t-\tau|^{2H}
\]
where $H$, $\bD$ are given by (\ref{061802}), (\ref{061803})
respectively.
\end{lemma}

{\bf Proof.} Thanks to the stationarity of the path $\bx_\ep(t)$
it is enough to prove the lemma for $\tau=0$. By the stationarity of
$\vV(s,\ep\bx(s))$ (\cite{port-stone}), we write 
\be
\label{1805}
\lim\limits_{\ep\da}\bE\left[\bx_\ep(t)\otimes\bx_\ep(t)\right]=
\lim\limits_{\ep\da}
\ep^2\int\limits_0^{\frac{t}{\ep^{2\delta}}}ds\int\limits_0^s
\bE\left[\vV(s',\ep\bx(s'))\otimes\vV(0,\bze)\right]ds'
\ee
which
equals
\be
\label{381}
2\sum\limits_{n=1}^N{\cal I}_n+{\cal R}_N
\ee
where
\[
{\cal
  I}_n=\ep^{n+1}\int\limits_0^{\frac{t}{\ep^{2\delta}}}ds
\int\limits_0^sds_1\cdots\int\limits_0^{s_{n-1}}
\bE\left[\bW_{n-1}(s_1,\cdots,s_n,\bze)\otimes\vV(0,\bze)\right]ds_n
\]
and
\[
\bW_0(s_1,\bx)=\vV(s_1,\bx)
\]
\[
\bW_{n}(s_1,\cdots,s_{n+1},\bx)=\vV(s_{n+1},\bx)\cdot\nabla\bW_{n-1}(s_1,\cdots,s_{n},\bx)\qquad\mbox{for
  }n=1,2,\cdots
\]
with the remainder term
\be
\label{5145}
{\cal
 R}_N=2\ep^{N+2}\int\limits_0^{\frac{t}{\ep^{2\delta}}}ds
\int\limits_0^sds_1\cdots\int\limits_0^{s_{N}}
\bE\left[\bW_{N}(s_1,\cdots,s_{N+1},\ep\bx(s_{N+1}))\otimes\vV(0,\bze)\right]ds_{N+1}.
\ee

{\em Estimates of ${\cal I}_n$.} Elementary calculations show that
\be
\label{I1}
\lim\limits_{\ep\da}{\cal I}_1=\bD t^{2H},\quad\hbox{for}\,\,\alpha+\beta>1.
\ee
Since $\vV$ is Gaussian  we have that
\[
\bE{\cal I}_n=\bze,\quad\hbox{for even $n$}.
\]
We now show that
\be
\label{b1}
\lim\limits_{\ep\da}\bE{\cal I}_n=\bze, \quad\hbox{for odd $n$}.
\ee

Set
\[
\bE_{s_{n+1}}W_{n-1,i}(s_1,\cdots,s_n,\bx):=\bE
\left[W_{n-1,i}(s_1,\cdots,s_n,\bx)\left|\right.{\cal
V}_{-\infty,s_{n+1}}\right].
\]
The $i,j$-th entry of the matrix ${\cal I}_n$ is given by
\be
\label{2211}
\ep^{n+1}\int\limits_0^{\frac{t}{\ep^{2\delta}}}ds
\int\limits_0^{s}ds_1\cdots\int\limits_0^{s_{n-1}}
\bE\left[\bE_0W_{n-1,i}(s_1,\cdots,s_n,\bze)V_j(0,\bze)\right]ds_n.
\ee
The conditional expectation in (\ref{2211}) can be expressed 
in terms of spectral measures of the velocity field. To do so we
introduce first the so-called {\em proper} functions of order $n,
 \bbsig:\{1,\cdots,n\}\rightarrow \{0,1\}$ that appear in the
statement of the next lemma. The proper function of order 1
is unique and is given
by $\bbsig(1)=0$. Any proper function, $\bbsig'$, of order
$n+1$ is generated
from a proper function $\bbsig$ of order $n$ as follows.
For some $p\leq n$,
\beq
\label{351}
\nonumber
\bbsig'(n+1)&:=&0\\
\nonumber
\bbsig'(k)&:=&\bbsig(k)\quad
\mbox{for}\quad k\leq n\mbox{ and }k\not=p\\
\bbsig'(p)&:=&1-\bbsig(p).
\eeq
In other words, each proper function $\bbsig$ of order $n$ generates $n$
different proper functions of order $n+1$. Thus, the total number of proper
functions of order $n$ is $(n-1)!$. In the sequel, we sometimes write $\sigma_k$
instead of $\sigma(k)$.

\begin{lemma}
\label{lemma6171}
Let $n\geq1$ and $s_1\geq s_2\geq \cdots \geq s_n\geq  s_{n+1}$,
$i\in\{1,\cdots,d\}$, $\bx\in R^d$. We have then that
\be
\label{1823}
\bE_{s_{n+1}}W_{n-1,i}(s_1,\cdots,s_n,\bx)
=
\ee
\[
\sum\int\cdots\int\varphi_{\bbi,\bbsig}^{(n)}(\bk_1,\cdots,\bk_n)
 \exp\{-\sum\limits_{m\in A_n({\cal
    F})}|\bk_m|^{2\bt}(s_n-s_{n+1})\}\times
\]
\[
P_{n-1}({\cal
    F})Q({\cal
    F})\prod\limits_{m\in A_n({\cal F})}
\widehat{V}_{i_{m},\si_{m}}(s_{n+1},\bx,d\bk_{m}),
\]
where $\vphi^{(n)}_{\bbi,\bbsig}$ are some functions with
 $\sup |\vphi^{(n)}_{\bbi,\bbsig}|\leq1$
\be
\label{a.4}
P_{n-1}({\cal
    F})=\prod\limits_{j=1}^{n-1}\left(\sum\limits_{m\in A_j({\cal
    F})}|\bk_m|\right)\exp\{-\sum\limits_{m\in A_j({\cal
    F})}|\bk_m|^{2\bt}(s_j-s_{j+1})\}
\ee
\be
\label{1833}
Q({\cal
    F})= \prod\limits_{\widehat{mm'}\in E({\cal
    F})}\bE\left[\widehat{V}_{i_{m},\si_{m}}(0,d\bk_{m})
\widehat{V}_{i_{m'},\si_{m'}}(0,d\bk_{m'})\right].
\ee
The  summation is over all multiindices $\bbi$ of length $n$,
whose first component equals $i$, all
 ${\cal F}\in \cG_s$ and all proper functions $\bbsig$ of order $n$.
\end{lemma}

Before proving Lemma~\ref{lemma6171}, we apply it to show  (\ref{b1}).
By (\ref{1823})
\[
\int\limits_0^{\frac{t}{\ep^{2\delta}}}ds
\int\limits_0^sds_1\cdots\int\limits_0^{s_{n-1}}
\bE_0W_{n-1,i}(s_1,\cdots,s_n,\bze)ds_n=
\]
\[
\sum\int\limits_0^{\frac{t}{\ep^{2\delta}}}ds
\int\cdots\int \tilde{\vphi}_{\bbi,\bbsig}^{(n)}
(\bk_1,\cdots,\bk_n)\exp\{-\sum\limits_{m\in A_n({\cal
    F})}|\bk_m|^{2\bt}s\} P_{n-1}'({\cal
    F})Q({\cal
    F})
\prod\limits_{m\in
    A_n({\cal
    F})}\widehat{V}_{i_m,\si_m}(0,d\bk_m)
\]
for $i=1,\cdots,d$. Here, adopting the convention $s_{n+1}:=0$, we set
\[
\widetilde{\varphi}_{\bbi,\bbsig}^{(n)}(\bk_1,\cdots,\bk_n):=
\]
\[
{{\int\limits_0^sds_1\cdots\int\limits_0^{s_{n-1}}
\varphi_{\bbi,\bbsig}^{(n)}(\bk_1,\cdots,\bk_n)\times\prod\limits_{j=1}^n
\exp\left\{-\sum\limits_{m\in A_j({\cal
    F})}|\bk_m|^{2\bt}(s_j-s_{j+1})\right\}ds_1\cdots ds_n}\over{
\int\limits_0^sds_1\cdots\int\limits_0^{s}\prod\limits_{j=1}^n
\exp\left\{-\sum\limits_{m\in A_j({\cal
    F})}|\bk_m|^{2\bt}s_j\right\}ds_1\cdots ds_n}}
\]
and
\[
P_{n-1}'({\cal
    F})=\prod\limits_{j=1}^{n-1}\left\{\left(\sum\limits_{m\in A_j({\cal
    F})}|\bk_{m}|\right)
\times{{1-\exp\{-\sum\limits_{m\in A_j({\cal
    F})}|\bk_m|^{2\bt}s\}}\over{\sum\limits_{m\in A_j({\cal
    F})}|\bk_{m}|^{2\bt}}}\right\}.
\]
It is elementary to check that, due to
$|\varphi_{\bbi,\bbsig}^{(n)}|\leq 1$,
\be
\label{1903}
|\widetilde{\varphi}_{\bbi,\bbsig}^{(n)}(\bk_1,\cdots,\bk_n)|\leq 1.
\ee

By Lemma \ref{lemma6171} the left hand side of (\ref{2211})
equals
\be
\label{2212}
2\ep^{n+1}
\sum\int\limits_0^{\frac{t}{\ep^{2\delta}}}ds\int\cdots\int
\frac{\tilde{\vphi}^{(n)}_{\bbi,\bbsig}(\bk_1,\cdots,\bk_n)}
{\sum\limits_{m\in A_n({\cal
    F})}|\bk_{m}|^{2\bt}} P_{n-1}'({\cal
    F})Q({\cal
    F})\bE\left[\prod\limits_{m\in
    A_n({\cal
    F})\cup\{n+1\}}\widehat{V}_{i_m,\si_m}(0,d\bk_m)\right].
\ee
Here the summation extends over all multiindices $\bbi=(i_1,\cdots,i_{n+1})$
 such that $i_1=i$,
$i_{n+1}=j$, all Feynman diagrams ${\cal F}\in \cG_s$ and all proper functions
$\bbsig$ of order $n$. Note that
\be
\label{22'}
\frac{1-e^{-\xi/\ep^{2\delta}}}{\xi}\leq
\frac{C}{\ep^{2\delta}+\xi}
\ee
for all positive $\ep,\xi$. Here and in the sequel $C$ stands
for a generic constant independent of $\ep$. $C$ in (\ref{22'}) is
also independent of $\xi>0$.
Thus, the absolute
value of (\ref{2212}) is bounded by
\be
\label{2213}
t\ep^{n+1-2\delta}
\sum\int\limits_0^{K}\cdots
\int\limits_0^{K}
\frac{p_{n-1,\ep}({\cal
    F})}
{\ep^{2\delta}+\sum\limits_{m\in A_n({\cal F})}k_m^{2\bt}}
\prod\limits_{\widehat{mm'}}\frac{\delta(k_m-k_{m'})
dk_mdk_{m'}}{k_m^{2\al-1}}
\ee
with
\[
p_{n-1,\ep}({\cal
    F}):=\prod\limits_{j=1}^{n-1}
    \frac{\sum\limits_{m\in A_j({\cal F})}k_m}
    {\ep^{2\delta}+\sum\limits_{m\in A_j({\cal F})}k_m^{2\bt}}.
\]

Using the fact that
\be
\label{012011}
\frac{\sum\limits_{m\in A_j({\cal
 F})}k_m}{\ep^{2\delta}+\sum\limits_{m\in A_j({\cal
 F})}k_m^{2\bt}}\leq
C\frac{\ep^{\frac{\delta}{\bt}}+k_{m_j}}{\ep^{2\delta}+k_{m_j}^{2\bt}},
\quad \forall m_j\in A_j(\cal F),
\ee
we have
\be
\label{224}
\frac{p_{n-1,\ep}({\cal F})}{\ep^{2\delta}+\sum\limits_{m\in A_n({\cal
F})}k_m^{2\bt}}\leq C\prod\limits_{j=1}^{n-1}
\frac{\ep^{\frac{\delta}{\bt}}+k_{m_j}}{\ep^{2\delta}+k_{m_j}^{2\bt}}
\frac{1}{\ep^{2\delta}+k_{m_n}^{2\bt}}, \quad \forall m_j\in A_j(\cal F).
\ee

 Bounds (\ref{224}) and (\ref{2213}) imply that
\be
\label{2252}
|{\cal I}_n|\leq Ct\ep^{n+1-2\delta}
\sum\limits_{\cF\in \cG_c}
\prod\limits_{\widehat{mm'}\in E({\cal F})}
\int\limits_0^{K}
\frac{\left(
\ep^{\frac{\delta}{\bt}}+k_{m}\right)^{q_m}}
{\left(\ep^{2\delta}+k_{m}^{2\bt}\right)^{q_m+\delta_{m,m_n}}}
\frac{dk_m}{k_m^{2\al-1}}
\ee
where $q_m$ are certain nonnegative exponents satisfying
\be
\label{5141}
\sum\limits_{\widehat{mm'}\in E({\cal F})} q_{m}=n-1
\ee
and $\delta_{m,m_n}=0$ if $m\not=m_n$ and $=1$ otherwise. 

%with $c_n$ denoting the cardinality of $A_n({\cal F})$. The summation above
%extends over all free vertices $m\in A_n(\cF)$.

The integrals appearing in the expression (\ref{2252}) are of the form
\be
\label{51410}
\int\limits_0^{K}
\frac{\left(\ep^{\frac{\delta}{\bt}}+k\right)^q}
{\left(\ep^{2\delta}+k^{2\bt}\right)^{q+r}}\frac{dk}{k^{2\al-1}}
\ee
for some $q\geq0$ and $r\in\{0,1\}$. They may diverge or remain bounded as
$\ep\da$ depending on the exponents $q,r$. If $q,r$ are such that the integral
diverges then $2\bt(q+r)+2\al>2+q$ and, consequently,
\[
\int\limits_0^{+\infty}
\frac{(1+k)^q}
{(1+k^{2\bt})^{q+r}}\times\frac{dk}{k^{2\al-1}}<+\infty.
\]
In either case, the integral (\ref{51410}) is bounded
from above by $C\ep^{n(2-\alpha-2\beta)/(\alpha+2\beta-1)}$
for $1<\alpha+2\beta$
after a change of variable $k'_j=k_j/\ep^{\frac{\delta}{\bt}}$
in case (\ref{51410}) diverges as $\ep\to 0$.
Therefore
\be
\label{22222}
|{\cal I}_n|\leq Ct\ep^{n+1-2\delta}
\ep^{\frac{n(2-\al-2\bt)}{\al+2\bt-1}-1}
%\sum
%\prod\limits_{\widehat{mm'}\in E({\cal F})}
%\int\limits_0^{K/\ep^{\frac{\delta}{\bt}}}
%\frac{(1+k_{m})^{q_m}}{(1+k_{m}^{2\bt})^{q_m+r_{m,m_n}}}
%\frac{dk_m}{k_m^{2\al-1}}\\
\leq
Ct\ep^{\frac{n-2}{\al+2\bt-1}}
\ee
which  vanish as $\ep\da$ for
$n\geq3$.

{\em Estimates of ${\cal R}_N$.} By (\ref{5145})
\[
{\cal
  R}_N=2\ep^{N+2}\int\limits_0^{\frac{t}{\ep^{2\delta}}}ds
\int\limits_0^sds_1\cdots\int\limits_0^{s_{N}}
\bE\left[\bE_{s_{N+1}}
\bW_{N}(s_1,\cdots,s_{N+1},\ep\bx(s_{N+1}))
\otimes\vV(0,\bze)\right]ds_{N+1}.
\]

By the Cauchy-Schwartz inequality we get that
\be
\label{1902}
|{\cal
  R}_N|^2\leq4t^2\ep^{4(1-\delta)+2N}\bE|\vV(0,\bze)|^2\times
\ee
\[
\max\limits_{0\leq s\leq t/\ep^{2\delta}}
\bE\left|
\mathop{\int\cdots\int}\limits_{
 s\geq s_1\geq\cdots\geq
    s_{N+1}\geq0}\bE_{s_{N+1}}\bW_{N}(s_1,\cdots,s_{N},s_{N+1},
\ep\bx(s_{N+1}))ds_1\cdots ds_{N+1}\right|^2.
\]
The stationarity of the Lagrangian velocity field implies that the maximum in
(\ref{1902}) is  equal to
\be
\label{3121}
\max\limits_{0\leq s\leq t/\ep^{2\delta}}
\bE\left|
\int\limits_0^{s}ds'
\mathop{\int\cdots\int}\limits_{
s'\geq s_1\geq\cdots\geq s_{N}\geq 0}
\bE_0\bW_{N}(s_1,\cdots,s_{N},0,\bze)ds_1\cdots ds_N\right|^2\leq
\ee
\[
C\max\limits_{0\leq s\leq t/\ep^{2\delta}}
\bE\left|
\int\limits_0^{s}ds'
\mathop{\int\cdots\int}\limits_{
s'\geq s_1\geq\cdots\geq s_{N}\geq 0}
\bE_0\nabla\bW_{N-1}(s_1,\cdots,s_{N-1},s_N,\bze)ds_1\cdots ds_N\right|^2\times
\bE\left|\vV(0,\bze)\right|^2.
\]
Here the hypercontractive property of 
the Gaussian measure is used. Subsequent applications of
Lemma \ref{lemma6171} to (\ref{3121})
yields the upper bound
\be
\label{1906}
C\frac{t^2}{\ep^{4\delta}}
\bE\left|\sum\limits_{\cF\in \cG_s,\sigma, \bbi}\int\cdots\int
\psi_{\bbi,\bbsig}(\bk_1,\cdots,\bk_{N}) P_N({\cal
 F})Q({\cal
 F})
\prod\limits_{m\in
A_{N}({\cal F})}\widehat{V}_{i_m,\si_m}(0,d\bk_m)
\right|^2
\ee
with some $|\psi_{\bbi,\bbsig}|\leq 1$. The summation above extends over all
Feynman diagrams ${\cal
F}\in\cG_s$, the relevant proper functions $\bbsig$  and multiindices $\bbi$.

Thus, we have
\be
\label{2222}
{\cal R}_N^2\leq Ct^4\ep^{2N+4(1-2\delta)}
\sum_{\cF,\cF'}\int\limits_0^{K}\cdots
\int\limits_0^{K}
p_{N,\ep}({\cal
    F})p_{N,\ep}({\cal
    F}')\prod\limits_{\widehat{mm'}}\frac{\delta(k_m-k_{m'})
dk_mdk_{m'}}{k_m^{2\al-1}}.
\ee
Here the summation extends over all possible completions of  with ${\cal
F}\in\cG_s(\{1,\cdots,N\})$, ${\cal F}'\in\cG_s(\{N+1,\cdots,2N\})$. The product
is over all edges of any completion of ${\cal F}\cup{\cal F}'$. Arguing as for
(\ref{2252}) we obtain that 
\be
\label{5146}
|{\cal
  R}_N|^2\leq Ct^4\ep^{2N+4(1-2\delta)}
\sum_{\cF,\cF'}
\prod\limits_{\widehat{mm'}}
\int\limits_0^{K}\left(\frac{\ep^{\frac{\delta}{\bt}}+k_{m}}{\ep^{2\delta}+k_{m}^{2\bt}}
\right)^{q_m}
\frac{dk_m}{k_m^{2\al-1}}
\ee
for some $q_m\geq0$
with
\be
\label{5147}
\sum\limits_{\widehat{mm'}}q_m=2N.
\ee
Moreover,
\be
\label{4271}
\label{2225}
|{\cal
  R}_N|^2\leq Ct^4\ep^{2N+4(1-2\delta)}
\ep^{-\frac{2N(2-\al-2\bt)}{\al+2\bt-1}}\leq
Ct^4\ep^{\frac{2N}{\al+2\bt-1}+4(1-2\delta)}.
\ee
which
  vanishes as $\ep\da$ for a sufficiently large $N$.
In conclusion,  we proved that the left hand side of (\ref{1805}) tends to $\bD
t^{2H}$ as $\ep\da$, provided that $\al+\bt>1$ (see (\ref{I1})).

By the hypercontractivity property of the $L^p$ norms over Gaussian measures 
we also know that for any $p\geq1$ and $T>0$
there exists a constant $C>0$
\be
\label{4281}
\bE|\bmx_\ep(t)-\bmx_\ep(s)|^p\leq C(t-s)^{2Hp}
\ee
for any $T\geq t\geq s\geq0$, $\ep>0$.

{\bf Proof of Lemma \ref{lemma6171}.} We prove the lemma by
 induction. The case $n=1$  is obvious by choosing
$\vphi^{(0)}_{i}\equiv 1$. Suppose that the
result holds for $n$. For the sake of convenience we assume
with no loss of any generality that $s_{n+2}=0$, then
\be
\label{2261}
\bE_0W_{n+1,i}(s_1,\cdots,s_{n+1},\bx)=
\bE_0\left\{\vV(s_{n+1},\bx)\cdot\nabla
  \bE_{s_{n+1}}W_{n,i}(s_1,\cdots,s_{n},\bx)\right\}.
\ee
By virtue of the inductive assumption we can represent
$\bE_{s_{n+1}}W_{n,i}$ using (\ref{1823}) and as a result
 (\ref{2261}) becomes
\be
\label{1224}
\sum\bE_0\left[\int\cdots\int\varphi_{\bbi,\bbsig}^{(n)}
(\bk_1,\cdots,\bk_n)\exp\{-\sum\limits_{m\in A_n({\cal
    F})}|\bk_m|^{2\bt}(s_n-s_{n+1})\}P_{n-1}({\cal
    F})Q({\cal
    F})\right.
\ee
\[
\left.\widehat{\vV}_0(s_{n+1},\bx,d\bk_{n+1})\cdot
\nabla\left\{\prod\limits_{m\in
    A_n({\cal
    F})}\widehat{V}_{i_m,\si_m}(s_{n+1},\bx,d\bk_m)\right\}\right].
\]

To calculate (\ref{1224}) we decompose each $\widehat{V}_{\si,i}(s,\bx,d\bk)$
as
\be
\label{a.6}
\widehat{V}_{\si,i}(s,\bx,d\bk)=\widehat{V}^0_{\si,i}(s,\bx,d\bk)
+\widehat{V}^1_{\si,i}(s,\bx,d\bk)
\ee
where
\be
\label{a.3}
\widehat{V}^0_{\si,i}(s,\bx,d\bk)=
e^{-|\bk|^{2\bt}(s-t)}\widehat{V}_{\si,i}(t,\bx,d\bk)
\ee
is the orthogonal projection of $\widehat{V}_{\si,i}$ on ${\cal V}_{-\infty,t}$.
Expression (\ref{1224}) becomes
\be
\label{1225}
\sum\bE_0\left[\int\cdots\int
    \vphi_{\bbi,\bbsig} ^{(n)}(\bk_1,\cdots,\bk_n)
\exp\{-\sum\limits_{m\in A_n({\cal F})}|\bk_m|^{2\bt}(s_n-s_{n+1})\}
P_{n-1}({\cal
    F})Q({\cal
    F}){\cal K}({\cal
    F}) \right]
\ee
with
\[
{\cal K}({\cal
    F}) := \sum\limits_{
\varrho=\{\varrho_{j}\}\atop j\in A_n(\cF)\cup\{n+1\}}
\hspace*{.005in}~\widehat{\vV}_0^{\varrho_{n+1}}(s,\bx,d\bk_{n+1})\cdot\nabla\left\{\prod\limits_{m\in
    A_n({\cal
    F})}\widehat{V}_{\si_m,i_m}^{\varrho_{m}}(s,\bx,d\bk_m)\right\}.
    \]
The term corresponding to $\varrho_j\equiv 1$
vanishes,
as is clear from the following calculation,
\be
\label{0311}
\bE\left\{\int\cdots\int \vphi_{\bbi,\bbsig} ^{(n)}
(\bk_1,\cdots,\bk_n)P_{n-1}({\cal
    F})Q({\cal
    F})
\right.
\ee
\[
\left.
\widehat{\vV}^{1}_0(s,\bx,d\bk_{n+1})
\cdot\nabla\left(\prod\limits_{m\in
    A_n({\cal
    F})}
\widehat{V}_{\si_m,i_m}^{1}(s,\bx,d\bk_m)\right)\right\}=
\]
\[\nabla\cdot\bE\left\{\int\cdots\int
\vphi_{\bbi,\bbsig} ^{(n)}(\bk_1,\cdots,\bk_n)P_{n-1}({\cal
  F})Q({\cal
  F})    \widehat{\vV}^{1}_0(s,\bx,d\bk_{n+1})\prod\limits_{m\in
    A_n({\cal
    F})}\widehat{V}_{\si_m,i_m}^{1}(s,\bx,d\bk_m)\right\}=0
\]
by homogeneity of the velocity field.
By (\ref{a.1})-(\ref{a.2})
\beq
&&\widehat{\vV}_0(s,\bx,d\bk_{n+1})\cdot\nabla\left\{\prod\limits_{m\in
A_n({\cal
F})}\widehat{V}_{\si_m,i_m}(s,\bx,d\bk_m)\right\}\label{a.7}\\
&=&\sum\limits_{m'\in
A_n({\cal
F})}\bk_{m'}\cdot\widehat{\vV}_{0}(s,\bx
,d\bk_{n+1})\prod\limits_{m\in
A_n({\cal
F})}\widehat{V}_{\si_m^{m'},i_m}(s,\bx,d\bk_m)
\nonumber
\eeq
where
\be
\label{46'}
\si_{m}^{m'}:=\left\{
\begin{array}{ll}
1-\si_{m'} &\mbox{if\quad}m'=m\\
\si_m &\mbox{otherwise.}
\end{array}\right.
\ee

By (\ref{a.3}), (\ref{a.6}), (\ref{a.7}) and the definition
(\ref{a.4}),  (\ref{1225}) further reduces to
\be
\label{061605}
\sum
\int\cdots\int\sum\limits_{i_{n+1}=1}^d
\sum\limits_{m'\in A_n({\cal F})}\sum\limits_{{\cal F}'}
\vphi_{\bbi,\bbsig} ^{(n)}
\frac{k_{m', i_{n+1}}}{\sum\limits_{m\in A_n({\cal F})}
|\bk_{m}|}
\ee
\[
\exp\{-\sum\limits_{m\in A({\cal
 F}')}|\bk_m|^{2\bt} s_{n+1}\}P_{n}({\cal
 F})Q({\cal
 F})
 \prod\limits_{m\in
  A({\cal F}')}\widehat{V}_{\si_m^{m'},i_m}(t,\bx,d\bk_m)
  \]
  \[
\left.\mathop{\prod\limits_{\widehat{pq}\in E({\cal
F}')}}\left[1-e^{-(|\bk_{p}|^{2\bt}+|\bk_{q}|^{2\bt})(s-t)}\right]\bE\left[
\widehat{V}_{\tilde{\si}_{p,m'},i_{p}}
(0,\bze,d\bk_{p})\widehat{V}_{\tilde{\si}_{q,m'
 },i_{q}}(0,\bze,d\bk_{q})\right]\right\}
 \]
with $\tilde{\sigma}_{1,m'}=0$ and $\tilde{\sigma}_{j+1,m'}=\sigma_j^{m'}$
and
  all incomplete Feynman diagrams ${\cal F}'$ based on the
  set $A_n({\cal F})\cup \{n+1\}$.

Lemma~2 follows with
\[
\vphi_{\bbi,\bbsig'_{m'}}^{(n+1)}(\bk_1,\cdots,\bk_{n+1}):=
\vphi_{\bbi,\bbsig} ^{(n)}(\bk_1,\cdots,\bk_{n})
\frac{k_{m',i_{n+1}}}{\sum\limits_{m'\in A_n({\cal F})}
  |\bk_{m'}|}
\mathop{\prod\limits_{\widehat{pq}\in E({\cal
    F}')}}\left[1-e^{-(|\bk_{p}|^{2\bt}+|\bk_{q}|^{2\bt})s_{n+1}}\right]
\]
\kropa

\section{Proof of weak convergence}
It is easy to see that the Gaussian processes
\be
\label{062502}
\by_\ep(t):=\ep\int\limits_0^{\frac{t}{\ep^{2q}}}\vV(s,\bze)ds~~t\geq0.
\ee
converge weakly to the fractional Brownian Motion $\bB_H(t)$, $t\geq0$ given by
(\ref{062701}). In addition we have
\[
\limsup\limits_{\ep\da}\bE|\by_\ep(t)|^p<+\infty
\]
for any $p\geq 1$, $t\geq0$.

We now prove that
\beq
\nonumber
&&\lim\limits_{\ep\da}\bE
\left\{[x_{\ep,i_1}(t_1)-x_{\ep,i_1}(t_2)]^{p_1}
\cdots[x_{\ep,i_M}(t_M)-x_{\ep,i_{M}}(t_{M+1})]^{p_M}\right\}\\
&=&\label{062901}
\bE\left\{[B_{H,i_1}(t_1)-B_{H,i_1}(t_2)]^{p_1}\cdots
[B_{H,i_M}(t_M)-B_{H,i_{M}}(t_{M+1})]^{p_M}\right\}
\eeq
which, in conjuction with the tightness, identifies
the fractional Brownian motion $B_H(t)$ as the limit.
Equation (\ref{062901}) is a consequence of
\be
\label{star}
\lim\limits_{\ep\da}\left|\bE\{[x_{\ep,i_1}(t_1)-x_{\ep,i_1}(t_2)]^{p_1}\cdots[x_{\ep,i_M}(t_M)-x_{\ep,i_{M}}(t_{M+1})]^{p_M}-\right.
\ee
\[
\left.[y_{\ep,i_1}(t_1)-y_{\ep,i_1}(t_2)]^{p_1}\cdots[y_{\ep,i_M}(t_M)-y_{\ep,i_{M}}(t_{M+1})]^{p_M}\}\right|=0
\]
with $\by_{\ep}(t)=(y_{\ep,1}(t),\cdots,y_{\ep,d}(t))$, which, in turn, follows
from the next lemma.
\begin{lemma}
\label{lemma6251}
For any positive integers $M$, $p_j$,  multiindices
$\bbi_j\in
\{1,\cdots,d\}^{p_j}$ with $j=1,\cdots,M$ 
we have
\be
\label{08175}
\lim\limits_{\ep\da}\left|\bE\left[Z_{\ep,\bbi_1}^{(p_1)}(t_2,t_1)
\cdots
Z_{\ep,\bbi_M}^{(p_M)}(0,t_M)-W_{\ep,\bbi_1}^{(p_1)}
(t_2,t_1)\cdots
    W_{\ep,\bbi_M}^{(p_M)}(0,t_{M})\right]\right|=0.
\ee
Here for any integer $N\geq1$,  multiindex $\bbi=(i_1,\cdots,i_N)\in
\{1\cdots,d\}^N$ and $t\geq s$  we define
\[
Z_{\ep,\bbi}^{(N)}(s,t):=\ep^N
\mathop{\int\cdots\int}\limits_{\triangle_N(s,t)}
\prod\limits_{p=1}^N V_{i_p}(s_p,\ep\bx(s_p))ds_1\cdots ds_N
\]
and
\[
W_{\ep,\bbi}^{(N)}(s,t):=\ep^N\mathop{\int\cdots\int}
\limits_{\triangle_N(s,t)}
\prod\limits_{p=1}^N V_{i_p}(s_p,\bze)ds_1\cdots ds_N,
\]
with $\triangle_N(s,t):=\{(s_1,\cdots,s_N):t/\ep^{2\delta}\geq s_1\geq
\cdots\geq s_N\geq s/\ep^{2\delta}\}$.
\end{lemma}

{\bf Proof.} To avoid  cumbersome expressions that may obscure the essence of
the proof we consider only the special case of $M=1$ and $t_1=t$, $t_2=0$. The
general case follows from exactly the same argument. We shall proceed with the
induction argument on $p_1=P$. The case when $P=1$ is trivial because the
stationarity of the relevant processes implies that the expression under the
limit in (\ref{08175}) vanishes. 
By (\ref{4281}) we know that
\[
\limsup\limits_{\ep\da}\bE|Z_{\ep,\bbi}^{(1)}(0,t)|^q<+\infty,
\quad \forall q\geq1.
\]
Suppose that (\ref{08175}) is true for $P\geq 2$ and that
\be
\label{5151}
\limsup\limits_{\ep\da}\bE|Z_{\ep,\bbi}^{(P-1)}(0,t)|^q<+\infty,
\quad \forall q\geq1.
\ee
Like (\ref{381}) we write 
\be
\label{382}
\bE Z_{\ep,\bbi}^{(P)}(0,t)=
\sum\limits_{n=0}^{N-1}
{\cal I}_{n}(0,t)+{\cal R}_{N}(0,t)
\ee
with
\be
\label{3811}
{\cal
I}_{n}(0,t):=\ep^{P+n+1}
\mathop{\int\cdots\int}\limits_{\Delta^{(n)}_P(0,t)}
\bE\left\{\bE_{s_2}W^{n}_{i_1}(\bbs_1^{(n)},\ep\bx(s_2))
\prod\limits_{p=2}^PV_{i_p}(s_p,\ep\bx(s_p))
\right\}d\bbs_1^{(n)}
ds_2\cdots ds_P
\ee
\be
\label{38111}
{\cal R}_{N}(0,t):=\ep^{P+N+1}\mathop{\int\cdots\int}
\limits_{\Delta^{(N)}_P(0,t)}
\bE\left\{\bE_{s_{1,N+1}}W^{N}_{i_1}(\bbs_1^{(N)},\ep\bx(s_{1,N+1}))
\prod\limits_{p=2}^PV_{i_p}(s_p,\ep\bx(s_p))\right\}d\bbs_1^{(N)}
ds_2\cdots ds_P.
\ee
Here
\[
\Delta^{(n)}_P(s,t):=\{(\bbs_1^{(n)},s_2,\cdots,s_P):
t/\ep^{2\delta}\geq \bbs_1^{(n)}
\geq s_2
\cdots\geq s_P\geq s/\ep^{2\delta}\}
\]
with $\bbs^{(n)}_1:=(s_{1,1},\cdots,s_{1,n+1})$. We write $t\geq \bbs_1\geq
s$, if $t\geq s_1\geq s_n\geq s$, where
$\bbs=(s_1,\cdots,s_n)$ is any ordered $n-$tuple in the sense that
$s_1\geq\cdots\geq s_n$.
\commentout{\[
\sum\limits_{n=0}^{N-1}
\ep^{n+P+1}\mathop{\int\cdots\int}\limits_{\Delta^{(n)}_P(0,t)}
\bE\left\{W_{n,i_1}(\bbs_1^{(n)},\ep\bx(s_2))
\prod\limits_{p=2}^PV_{i_p}(s_p,\ep\bx(s_p))d\bbs_1^{(n)}
ds_2\cdots ds_P+
\]
\[
\ep^{P+N+1}\mathop{\int\cdots\int}
\limits_{\Delta^{(n)}_P(0,t)}
\bE\left\{W_{N,i_1}(\bbs_1^{(N)},\ep\bx(s_{1,N+1}))
\prod\limits_{p=2}^PV_{i_p}(s_p,\ep\bx(s_p))
\right\}d\bbs_1^{(N)}
ds_2\cdots ds_P
\]
is most conveniently
stated in terms of the so-called {\em ordered multituples}. We say that an
$m$-tuple $\bbs:=(s_1,\cdots,s_m)\in R^m$ is {\em ordered} if $s_1\geq\cdots\geq
s_m$. For any integrals $m,q$ and two ordered multituples
$\bbs:=(s_1,\cdots,s_m)$ and $\bbt:=(t_1,\cdots,t_q)$ we write that $\bbs\geq
\bbt$ if $\min s_i\geq \max t_j$.
(\ref{382}) can now be rewritten in terms of conditional expectations as follows
}

The argument of the proof of Lemma \ref{lemma6171}  and
(\ref{5151}) imply that
\[
\lim\limits_{\ep\da}{\cal I}_{n}(0,t)=0,\quad n\geq 1
\]
and
\[
\lim\limits_{\ep\da}{\cal R}_{N}(0,t)=0
\]
for $N$ sufficiently large. Thus
$\bE
Z_{\ep,\bbi}^{(P)}(0,t)$ has the same limit as the term
\be
\label{5152}
{\cal I}_{0}(0,t):=\ep^{P+1}\mathop{\int\cdots\int}
\limits_{\Delta^{(0)}_P(0,t)}
\bE\left\{V_{i_1}(s_1,\ep\bx(s_2))
\prod\limits_{p=2}^PV_{i_p}(s_p,\ep\bx(s_p))
\right\}d s_1 \cdots
ds_P.
\ee
For (\ref{5152}) we use  a generalization of the argument
of the proof of Lemma \ref{lemma6171}. Let us introduce some additional
notation. For any multiindex $\bbi=(i_1,\cdots,i_p)$ and $p\geq1$ we define
$W^{p,n}_{\bbi}$ by induction as follows. We set
\[
W^{p,0}_{i_1,\cdots,i_p}(s_1,\cdots,s_p,\bx):= V_{i_1}(s_1,\bx)\cdots
V_{i_p}(s_p,\bx)-\bE\{ V_{i_1}(s_1,\bx)\cdots V_{i_p}(s_p,\bx)\}
\]
and 
\[
W^{p,n+1}_{i_1,\cdots,i_p}(s_1,\cdots,s_{p-1},
\bbs_p^{(n+1)},\bx):=\nabla W^{p,n}_{i_1,\cdots,i_p}(s_1,\cdots,
\bbs_p^{(n)},\bx)\cdot\vV(s_{p,n+2},\bx)
\]
 for any ordered $(n+1)$-tuple
 $\bbs_p^{(n)}=(s_{p,1},\cdots,s_{p,n+1})\leq s_{p-1}$ 
and $(n+2)$-tuple $\bbs_p^{(n+1)}=(s_{p,1},\cdots, s_{p,n+1},
s_{p,n+2})$. Expanding the left hand side of (\ref{5152}) like
(\ref{381}) we obtain that
\[
{\cal I}_{0}(0,t)=
\]
\[
\ep^{P+1}\mathop{\int\cdots\int}
\limits_{\Delta^{(0)}_P(0,t)}
\bE\left\{V_{i_1}(s_1,\ep\bx(s_2))
V_{i_2}(s_2,\ep\bx(s_2))\right\}\bE\left\{\prod\limits_{p=3}^P
V_{i_p}(s_p,\ep\bx(s_p))\right\}d s_1 ds_2\cdots ds_P+
\]
\[
\sum\limits_{n=0}^{N-1} {\cal I}_{1,n}(0,t)+{\cal R}_{1,N}(0,t)
\]
where
\be
\label{515111}
{\cal
I}_{1,n}(0,t):=
\ee
\[
\ep^{P+n+1}\mathop{\int\cdots\int}
\limits_{\Delta^{(1,n)}_P(0,t)}
\bE\left\{\bE_{s_3}W^{2,n}_{i_1,i_2}(s_1,\bbs_2^{(n)},\ep\bx(s_3))
\prod\limits_{p=3}^PV_{i_p}(s_p,\ep\bx(s_p))
\right\}ds_1d\bbs_2^{(n)}
ds_2\cdots ds_P
\]
\be
\label{5151111}
{\cal R}_{1,N}(0,t):=
\ee
\[
\ep^{P+N+1}\mathop{\int\cdots\int}
\limits_{\Delta^{(1,N)}_P(0,t)}
\bE\left\{\bE_{s_{2,N+1}}W^{2,N}_{i_1,i_2}(s_1,\bbs^{(N)},
\ep\bx(s_{2,N+1}))
\prod\limits_{p=3}^PV_{i_p}(s_p,\ep\bx(s_p))
\right\} ds_1 d\bbs_2^{(N)} ds_3\cdots ds_P,
\]
\[
\Delta^{(1,n)}_P(0,t):=\{(s_1,\bbs_2^{(n)},s_3,\cdots,s_P):t/\ep^{2\delta}\geq s_1
\geq\bbs_2^{(n)}\geq\cdots\geq s_P\geq0\}.
\]
We represent the conditional expectations appearing in (\ref{515111}) and
(\ref{5151111}) using a generalization (Lemma \ref{lemma6291}) of Lemma
\ref{lemma6171}.

To formulate it we need a generalized notion of a proper function, which we call
a $p$-proper function. Let $p$ be a positive integer. The $p$-proper function of
order 1 is unique and is given by $\bbsig(i)=0$, $i=1,\cdots,p$. Any $p$-proper
function, $\bbsig'$, of order $n+1$ is generated from a $p$-proper function
$\bbsig$ of order $n$ as follows. For some $q\leq p+ n$,
\beq
\label{35111}
\nonumber
\bbsig'(p+n+1)&:=&0\\
\nonumber
\bbsig'(k)&:=&\bbsig(k)\quad
\mbox{for}\quad k\leq n+p\mbox{ and }k\not=q\\
\bbsig'(q)&:=&1-\bbsig(q).
\eeq

We also distinguish a special class of Feynman diagrams $\cG^p_s(B)$
: a diagram $\cF$ of order $n+p$ belongs
to $\cG^p_s(B)$ if $A_k(\cF)$ is not empty for all $k=p,...,n+p.$

\begin{lemma}
\label{lemma6291}
For any positive integer $p$, $s_1\geq
\cdots\geq s_{p-1}\geq\bbs_p^{(n-1)}\geq s$, a multiindex
 $\bbi=(i_1,\cdots,i_p)\in\{1,\cdots,d\}^p$ we have
\be
\label{18231}
 \bE_sW^{p,n-1}_{\bbi}(s_1,\cdots,s_{p-1},\bbs^{(n-1)}_p,\bx)=\sum
\int\cdots\int\varphi_{\bbj,\bbsig}^{(p,n)}(\bk_1,\cdots,\bk_{p+n})
\ee
\[
\exp\{-\sum\limits_{m\in A_{n+p}({\cal
    F})}|\bk_{m}|^{2\bt}(s_{p,n}-s)\}
P_{p,n-1}({\cal F}) Q({\cal F})
 \prod\limits_{m\in A_{n+p}({\cal
    F}) }\widehat{V}_{i_{m},\si_{m}}(s,\bx,d\bk_{m}),
\]
where  $\varphi^{(p,n)}_{\bbj,\bbsig}$  are  functions  satisfying $
|\varphi^{(p,n)}_{\bbj,\bbsig}|\leq 1$ and
\[
P_{p,n}({\cal
    F})=\prod\limits_{j=p}^{n+p-1}
\left(\sum\limits_{m\in A_{j}({\cal
    F})}|\bk_{m}|\right) \exp\{-\sum\limits_{m\in
A_{j}({\cal
    F})}|\bk_{m}|^{2\bt}(s_{p,j-p}-s_{p,j-p+1})\},
\]
\be
\label{18331}
Q({\cal
    F})=\prod\limits_{\widehat{mm'}\in E({\cal
    F})}
\bE\left[\widehat{V}_{i_{m},\si_{m}}(0,d\bk_{m})
\widehat{V}_{i_{m'},\si_{m'}}(0,d\bk_{m'})\right].
\ee
The summation is over all multiindices $\bbj=(j_1,\cdots,j_{n+p})$, such that
$\bbj_{|p}=\bbi$,
 all ${\cal F}\in \cG^p_s$ and
 all $p$-proper functions $\bbsig$ of order $n$.
Here by a convention $s_{p,0}:=s_{p-1}$.
\end{lemma}
The proof of Lemma \ref{lemma6291} is exactly the same as that of Lemma
\ref{lemma6171} and is omitted.

Continue the proof of Lemma~3 using Lemma~4 
we have that ${\cal I}_{0}(0,t)$ is asymptotically equal to $\bE
Z_{\ep,\bbi}^{(P)}(0,t)$ and
\[
\ep^{P+1}\mathop{\int\cdots\int}
\limits_{\Delta_P(0,t)}
\bE\left\{V_{i_1}(s_1,\ep\bx(s_3))V_{i_2}(s_2,\ep\bx(s_3))
\prod\limits_{p=3}^PV_{i_p}(s_p,\ep\bx(s_p))\right\}ds_1
ds_2\cdots ds_P
\]
is asymptotically equal to $\bE Z_{\ep,\bbi}^{(P)}(0,t)$, as $\ep\da$. Repeating
the above argument $p$-times we obtain (\ref{08175}). Finally 
the hypercontractivity properties of the $L^p$ norms over Gaussian measure space
imply that (\ref{5151}) holds with $P-1$ replaced by $P$\kropa

\end{document}